\documentclass[12 pt,leqno]{amsart}
\usepackage{amsmath, amsthm, amscd, amsfonts, amssymb, graphicx, color}
\usepackage[bookmarksnumbered, plainpages]{hyperref}
\usepackage[margin=1in,headsep=.2in]{geometry}
\newtheorem{thm}{Theorem}[section]

\newtheorem{lem}[thm]{Lemma}
\newtheorem{prop}[thm]{Proposition}

\numberwithin{equation}{section}

\begin{document}                  

\title{Zeta functions and subgroup growth in $P2/m$}

\author{Hermina ALAJBEGOVI\'C \and  Muharem AVDISPAHI\'C }



\maketitle

\begin{abstract}
By means of zeta and normal zeta functions of space groups, we determine the number of subgroups, resp. normal 
subgroups, of the tenth crystallographic group for any given index. This enables us to draw conclusions on the
subgroup growth and the degree of this group.\\
\\
{\bf Mathematics Subject Classification (2010) 11M41, 20H15}\\
\end{abstract}

\section{ Introduction}

\mbox{}

 The zeta function of a group
$G$ is defined as ${\zeta _G}(s) = \sum\limits_{n \in \mathbb {N}}{{a_n}(G){n^{ - s}}} $, where ${a_n}(G)$ denotes the number of
subgroups of index $n$ in $G$. Analogously, the normal zeta
function of a group $G$ is given by $\zeta _G^ \triangleleft (s) =
\sum\limits_{n \in \mathbb {N}} {{c_n}(G){n^{ - s}}} $, where ${c_n}(G)$ is
the number of normal subgroups of index $n$ in $G$. These functions provide a useful tool for studying the relationship
between the asymptotic behavior of the sequences ${a_n(G)}$, resp. ${c_n(G)}$, and the structure of $G$. 
The concepts of the zeta and normal zeta function were applied to nilpotent groups by Smith \cite {S}, and
Grunewald, Segal and Smith \cite {GSS}. Building upon our previous results related to the space groups with the 
point group isomorphic
to the cyclic group of order 2 (see \cite{AA}), we derive explicit expressions for the zeta and normal zeta function 
of $P2/m$ in Sections 2. and 3., and determine the exact 
number of its subgroups and normal subgroups of finite index, in Section 4. In the final section, we turn our attention
to the subgroup growth. 

The group $P2/m$ is the tenth group in the International Tables for Crystallography \cite {H}. It contains translations, reflections and 
diad rotations. The translations form a normal abelian subgroup $T$ of rank 3 - the translation subgroup of $P2/m$ or 
the Bravais lattice. The point group of $P2/m$, i.e., its quotient by  the translation subgroup $T$, is a finite 
group isomorphic to the direct product of two cyclic groups of order two (Klein 4-group).

A minimal set of generators of $P2/m$ and the algebraic relations the generators satisfy are as follows (cf.\cite {L})
 \[G = P2/m = \left\langle {x,y,z,r,m\left| {\begin{array}{*{20}{c}}
  {\left[ {x,y} \right],\left[ {x,z} \right],\left[ {y,z} \right],{r^2},{m^2},{{(mr)}^2},{x^m} = x,} \\ 
  {{y^m} = {y^{ - 1}},{z^m} = z,{x^r} = {x^{ - 1}},{y^r} = y,{z^r} = {z^{ - 1}}} 
\end{array}} \right.} \right\rangle .\]

The subgroups ${G_{{2_1}}} = \left\langle {x,y,z,m} \right\rangle $, ${G_{{2_2}}} = \left\langle {x,y,z,r} \right\rangle $,
${G_{{2_3}}} = \left\langle {x,y,z,mr} \right\rangle $ of the group $P2/m$ are isomorphic to space groups $Pm$, $P2$ and $P{\bar1}$, respectively, 
while the subgroup ${G_3}= T = \left\langle {x,y,z} \right\rangle $ is isomorphic to $P1$, i.e., to ${\mathbb{Z}^3}$.
Therefore, the part of knowledge about the zeta and normal zeta functions of these groups,
summarized in the next theorem, will be useful for our present purpose.

For a sake of bravity, we denote the translates of the
Riemann zeta function by:\\ ${\zeta_k(s)=\zeta(s-k)}$, i.e., ${\zeta_2(s)=\zeta(s-2)}$.

Recall that ${\zeta _k}(s) = \sum\limits_{n \in \mathbb {N}}{{{n^{ - s+k}}}}$ converges absolutely for $Re(s)>k+1$ and
has a meromorphic extension to the whole complex plane with a simple pole at $s=k+1$.

\begin{thm} {\rm (see \cite{AA})} Zeta and normal zeta functions of groups $P{\bar1}$, $P2$ and $Pm$ read as follows 

$\zeta_{P{\bar 1}}(s)= {\zeta_1(s)\zeta_2(s)\zeta_3(s) + {2^{ - s}}\zeta(s)\zeta_1(s)\zeta_2(s)}$ 

$\zeta_ {P2}(s)= {(1 + {2^{ - s + 3}})\zeta (s)\zeta_1(s)\zeta_2(s)} $

$\zeta_ {Pm}(s)={(1 + 9 \cdot {2^{ - s}} + 6 \cdot {2^{ - 2s}})\zeta
(s)\zeta_(s)\zeta_1(s) + {2^{ - s}}\zeta (s)\zeta_1(s)\zeta_2(s)} $

$\zeta_{P{\bar 1}}^
\triangleleft (s)={1 + 14 \cdot {2^{ - s}} + 28 \cdot
{2^{ - 2s}} + 8 \cdot {2^{ - 3s}} + {2^{ -s}}\zeta(s)\zeta_1(s)\zeta_2(s)} $
 
$\zeta_{P2}^
\triangleleft (s)={(1 + 13 \cdot {2^{ - s}} + 22 \cdot {2^{ - 2s}} + 4 \cdot {2^{ -
3s}}) \cdot \zeta (s) + (3 \cdot {2^{ -2s}} + {2^{ - s}})\zeta (s)\zeta (s)\zeta_1(s)}$

$\zeta_{Pm}^\triangleleft (s) ={(1 + 11 \cdot {2^{ - s}} + 12 \cdot {2^{ - 2s}})\zeta
(s)\zeta_1(s) + {2^{ - s}}(1 + 3 \cdot {2^{ -
s}})\zeta (s)\zeta (s)\zeta_1(s)}$.

\end{thm}

Due to group isomorphisms mentioned above, the latter explicit expressions are the building blocks in forming the 
zeta and normal zeta function of $P2/m$.\\

\section{Zeta function of $P2/m$}
\mbox{}
\begin{thm}  The zeta  function of the space groups $P2/m$ is given by:

$\begin{array}{l}
{\zeta _{P2/m}}(s) = (1 + 20 \cdot {2^{ - s}} + 36 \cdot {2^{ - 2s}})\zeta _1^2(s){\zeta _2}(s) + {2^{ - s}} \cdot (1 + 9 \cdot {2^{ - s}} +\\+ 6 \cdot {2^{ - 2s}})\zeta (s)\zeta _1^2(s) +  {2^{ - s}}(1 + 8 \cdot {2^{ - s}})\zeta (s){\zeta _1}(s){\zeta _2}(s) + {2^{ - s}} \cdot {\zeta _1}(s){\zeta _2}(s){\zeta _3}(s).
\end{array}$

\end {thm}

{\bf Proof.} The proof proceeds in five steps. First, we count only those subgroups of  ${G_1} = \left\langle G \right\rangle$
that are not contained in  ${G_{{2_1}}} $, ${G_{{2_2}}} $, ${G_{{2_3}}}$, ${G_3}$. Then, we count those subgroups of 
${G_{{2_1}}}$ that are not contained in  ${G_{{2_2}}} $, ${G_{{2_3}}}$, ${G_3}$. The same procedure applies to 
${G_{{2_2}}} $ and ${G_{{2_3}}}$. This way, we avoid over-counting of subgroups of a finite index.

Now, any subgroup of  ${G_1}$ has the form ${H_1} = \left\langle {m{x^a}{y^b}{z^c},r{x^d}{y^e}{z^f},{x^g}{y^h}{z^i},{y^j}{z^k},{z^l}} \right\rangle $,
where $a, b, c, d, e, f, g, h, i, j, k$ and $l$ are integers. To avoid over - counting, we require that 
$0 \leqslant a,d < g;0 \leqslant b,e,h < j;0 \leqslant c,f,i,k < l$ \cite {M}. The index of this subgroup is $gjl$.
Note that we cannot allow $g$, $j$  or $l$ to be 0 as this would give a subgroup of infinite index in $G$. The restrictions on those possible values
are represented in the following tableau

$\left( {\begin{array}{*{20}{c}}
  1&0&a&b&c \\ 
  0&1&d&e&h \\ 
  0&0&g&h&i \\ 
  0&0&0&j&k \\ 
  0&0&0&0&l 
\end{array}} \right)$.\\

Reading
down the columns, this tableau quickly sums up the information we have just derived about $H_1$. 
 If $H_1$ is a subgroup of $G$, then according to the second isomorphism theorem, ${{H_1} \cap T}$ has to be normal subgroup in $H_1$ and ${{{H_1}} \mathord{\left/
 {\vphantom {{{H_1}} {{H_1} \cap T}}} \right.
 \kern-\nulldelimiterspace} {{H_1} \cap T}} \cong {{{H_1}T} \mathord{\left/
 {\vphantom {{{H_1}T} T}} \right.
 \kern-\nulldelimiterspace} T}$. For ${{H_1} \cap T}$ to be a normal subgroup in $H_1$, we 
must have ${u^{ - 1}}({H_1} \cap T)u \in {H_1} \cap T$ for $\forall u \in {H_1}$.
To verify this, it is sufficient to take the generators of $H_1$ and ${{H_1} \cap T}$. Let us take  
$u = m{x^a}{y^b}{z^c} \in {H_1}$ and ${x^g}{y^h}{z^i} \in {H_1} \cap T$. Now, we have:
${\left( {m{x^a}{y^b}{z^c}} \right)^{ - 1}}({x^g}{y^h}{z^i})\left( {m{x^a}{y^b}{z^c}} \right) = {z^{ - c}}{y^{ - b}}{x^{ - a}}{m^{ - 1}}{x^g}{y^h}{z^i}m{x^a}{y^b}{z^c}\\
= {z^{ - c}}{y^{ - b}}{x^{ - a}} ( {x^g}{y^{ - h}}{z^i}) {x^a}{y^b}{z^c} = {x^g}{y^{ - h}}{z^i}$.\\
Hereof, ${\left( {m{x^a}{y^b}{z^c}} \right)^{ - 1}}({x^g}{y^h}{z^i})\left( {m{x^a}{y^b}{z^c}} \right)\in {{H_1} \cap T}$, if  ${x^g}{y^{ - h}}{z^i}\in {{H_1} \cap T}$.

Repeating the process for the remaining generators, we get another condition  ${y^{- j}}{z^{k}}\in {{H_1} \cap T}$.

Since  ${{{H_1}} \mathord{\left/
 {\vphantom {{{H_1}} {({H_1} \cap T)}}} \right. \kern-\nulldelimiterspace} {({H_1} \cap T)}}$ is isomorphic to a 
subgroup of the Klein group, we see that
 ${(m{x^a}{y^b}{z^c})^2}$ ,${(r{x^d}{y^e}{z^f})^2}$, $ {( m{x^a}{y^b}{z^c}r{x^d}{y^e}{z^f})^2} \in {H_1} \cap T$.
 Using the relations between elements in the group, we conclude that the condition 
${(m{x^a}{y^b}{z^c})^2}$ $\in {H_1}\cap T$   is equivalent to ${x^{2a}}{z^{2c}}$  $\in {H_1} \cap T$. Indeed,

${\left( {m{x^a}{y^b}{z^c}} \right)^2} = m{x^a}{y^b}{z^c}m{x^a}{y^b}{z^c} = mm{x^a}{y^{ - b}}{z^c}{x^a}{y^b}{z^c} 
= {x^{2a}}{z^{2c}} \in {H_1} \cap T$.

The remaining conditions lead to another requirement ${y^{2e}}$ $ \in {H_1} \cap T$.

So, we end up with the following conditions ${x^g}{y^{ - h}}{z^i}$, ${x^{2a}}{z^{2c}}$, ${y^{ - j}}{z^k}$, ${y^{2e}}$ $ \in {H_1} \cap T$. 
If ${x^g}{y^{ - h}}{z^i}$ lies in ${H_1} \cap T$, then there exist integer numbers 
${\alpha _1}, {\beta_1}, {\gamma_1}$  such that  ${x^g}{y^{ - h}}{z^i}=({x^g}{y^h}{z^i})^{\alpha_1} ({y^j}{z^k})^{\beta_1}({z^l})^{\gamma_1}$.
Thus, we get the following system of equations \\

$C_1 = \left\{ \begin{gathered}
  g = g{\alpha _1}, - h = h{\alpha _1} + j{\beta _1},i = i{\alpha _1} + k{\beta _1} + l{\gamma _1}, \hfill \\
  2a = g{\alpha _2},0 = h{\alpha _2} + j{\beta _2},2c = i{\alpha _2} + k{\beta _2} + l{\gamma _2}, \hfill \\
  0 = g{\alpha _3}, - j = h{\alpha _3} + j{\beta _3},k = i{\alpha _3} + k{\beta _3} + l{\gamma _3}, \hfill \\
  0 = g{\alpha _4},2e = k{\alpha _4} + j{\beta _4},0 = i{\alpha _4} + k{\beta _4} + l{\gamma _4} \hfill \\ 
\end{gathered}  \right\}.$\\

By taking into account the conditions  $0 \leqslant a,d < g;0 \leqslant b,e,h < j;0 \leqslant c,f,i,k < l$,  this system can be reduced to\\

$C'_1 = \left\{ \begin{gathered}
   - 2h = j{\beta _1},0 = k{\beta _1} + l{\gamma _1}, \hfill \\
  2a = g{\alpha _2},0 = h{\alpha _2} + j{\beta _2},2c = i{\alpha _2} + k{\beta _2} + l{\gamma _2}, \hfill \\
  2k = l{\gamma _3}, \hfill \\
  2e = j{\beta _4},0 = k{\beta _4} + l{\gamma _4} \hfill \\ 
\end{gathered}  \right\}.$\\

To solve the system, we distinguish eight cases depending on the parity of each of the numbers $ g, j, l $.
We keep in mind that $0 \leqslant a,d < g;0 \leqslant b,e,h < j;0 \leqslant c,f,i,k < l$. So, if $ g, j, l $
  are odd numbers, we see that  $a$ has to be 0. Hence ${\alpha _2} = {\beta _2} = 0$. Since $l$ is odd, it follows
 that $c = 0$. Similarly, we get   $k = 0$  and  $e = 0$. From  ${\alpha _2} = 0$, it follows that there exist
  $l$ choices for $i$.  There are no additional restrictions on $b, d, f$. 
Thus, the contribution to the zeta function of group $P2/m$ coming from this case is: $\sum\limits_{g,j,l \in {\Bbb N}'} {{g^{ - s}}} {j^{ - s}}{l^{ - s}} \cdot g \cdot j \cdot {l^2}$.\\
In other seven cases, we get the contributions:\\\\
 $4 \cdot \sum\limits_{j,g \in \mathbb{N}',l \in 2\mathbb{N}} {{g^{ - s}}} {j^{ - s}}{l^{ - s}} \cdot g \cdot j \cdot {l^2}$,\hskip 3mm if $ g, j$ are odd and $l$ is even; \\
 $4 \cdot \sum\limits_{l,g \in \mathbb{N}',j \in 2\mathbb{N}} {{g^{ - s}}} {j^{ - s}}{l^{ - s}} \cdot g \cdot j \cdot {l^2}$,\hskip 3mm if $ g, l$ are odd and $j$ is even;\\ 
 $2 \cdot \sum\limits_{l,j \in \mathbb{N}',g \in 2\mathbb{N}} {{g^{ - s}}} {j^{ - s}}{l^{ - s}} \cdot g \cdot j \cdot {l^2}$,\hskip 3mm if $ l, j$ are odd  and $g$ is even;\\
 $6 \cdot \sum\limits_{l \in \mathbb{N}',g,j \in 2\mathbb{N}} {{g^{ - s}}} {j^{ - s}}{l^{ - s}} \cdot g \cdot j \cdot {l^2}$,\hskip 3mm if $ g, j$ are even  and $l$ is odd;\\
 $6 \cdot \sum\limits_{j \in \mathbb{N}',g,l \in 2\mathbb{N}} {{g^{ - s}}} {j^{ - s}}{l^{ - s}} \cdot g \cdot j \cdot {l^2}$,\hskip 3mm if $ g, l$ are even and $j$ is odd;\\
 $10 \cdot \sum\limits_{g \in \mathbb{N}',j,l \in 2\mathbb{N}} {{g^{ - s}}} {j^{ - s}}{l^{ - s}} \cdot g \cdot j\cdot {l^2}$,\hskip 2mm if $ l, j$ are even and $g$ is odd;\\
 $13 \cdot \sum\limits_{g,j,l \in 2\mathbb{N}} {{g^{ - s}}} {j^{ - s}}{l^{ - s}} \cdot g \cdot j \cdot {l^2}$,\hskip 6mm if $ i, j, g$ are even.\\

Adding the above contributions, we see that the total share in the zeta function of $P2/m$ coming from subgroups of the 
form $ H_1 $ is:  $(1 + 20 \cdot {2^{ - s}} + 36 \cdot {2^{ - 2s}})\zeta_2 (s)\zeta_1 (s)\zeta_1 (s)$.

If $H_2$ is a subgroup of  ${G_{{2_1}}}$, then $\left| {G:{H_2}} \right| = \left| {G:{G_{{2_1}}}} \right| \cdot \left| {{G_{{2_1}}}:{H_2}} \right| = 2 \cdot \left| {{G_{{2_1}}}:{H_2}} \right|$. 
Taking only those subgroups of ${G_{{2_1}}}$ that are not contained in  ${G_{{2_2}}} $, ${G_{{2_3}}}$, ${G_3} $ and
making use of the respective part of Theorem 1.1, we derive the following share in the zeta function coming from 
subgroups of the form $ H_2 $: $ 2^{ - s}(1 + 9 \cdot {2^{ - s}} +6 \cdot {2^{ - 2s}})\zeta _1^2(s){\zeta}(s)$.

Now, let $H_3$ be a subgroup of  ${G_{{2_2}}}$. Then $\left| {G:{H_3}} \right| = \left| {G:{G_{{2_2}}}} \right| \cdot \left| {{G_{{2_2}}}:{H_3}} \right| = 2 \cdot \left| {{G_{{2_2}}}:{H_3}} \right|$.
In view of Theorem 1.1., those subgroups of ${G_{{2_2}}}$ that are not contained in ${G_{{2_3}}}$ and ${G_3} $, yield the share: ${2^{ - s}} \cdot (1 + 7 \cdot {2^{ - s}})\zeta (s)\zeta_1 (s)\zeta_2 (s)$.

For a subgroup $H_4$ of the group ${G_{{2_3}}}$, we have  $\left| {G:{H_4}} \right| = \left| {G:{G_{{2_3}}}} \right| \cdot \left| {{G_{{2_3}}}:{H_4}} \right| = 2 \cdot \left| {{G_{{2_3}}}:{H_4}} \right|$.
Now, the subgroups of ${G_{{2_3}}}$ that are not contained in ${G_3} $, combined with the information from Theorem 1.1., imply the share: ${2^{ - s}} \cdot \zeta_1 (s)\zeta_2 (s)\zeta_3 (s)$.

Finally, we still have to consider the subgroups of the translation subgroup 
$T={G_3} = \left\langle {x,y,z} \right\rangle $.  If $H_5$ is a subgroup of  ${G_3}$, then 
$\left| {G:{H_5}} \right| = \left| {G:{G_3}} \right| \cdot \left| {{G_3}:{H_5}} \right| = 4 \cdot \left| {{G_3}:{H_5}} \right|$.
 Since the zeta function of $T \cong {\mathbb{Z}^3}$ is $\zeta(s)\zeta _1(s)\zeta_2(s)$,  
we get the share: $2^{ -2 s}\zeta (s){\zeta _1}(s){\zeta_2}(s)$.

Combining all above contributions stemming from subgroups $H_1$, $H_2$, $H_3$, $H_4$ and $H_5$, we get the zeta 
function of $P2/m$ as stated in the Theorem.\\

\section  {\bf Normal zeta function of $P2/m$}

\mbox{}
\begin {thm}  The normal zeta function of $P2/m$ is given by:\\
$\begin{array}{l}
  \zeta _{{P_{2/m}}}^ \triangleleft (s) =  1 + 29 \cdot {2^{ - s}} + 126 \cdot {4^{ - s}} + 92 \cdot {8^{ - s}} + 8 \cdot {16^{ - s}} + {2^{ - s}}(1 + 13 \cdot {2^{ - s}} + 22 \cdot {2^{ - 2s}} + \\+4 \cdot {2^{ - 3s}})\zeta (s) + {2^{ - s}}(1 + 11 \cdot {2^{ - s}} + 12 \cdot {2^{ - 2s}})\zeta (s){\zeta_1} (s) + {2^{ - 3s}}(3 + {2^s}){\zeta ^2}(s){\zeta_1} (s). \\ 
\end{array} $
\end {thm}

{\bf Proof.} We apply a similar procedure as in the case of Theorem 2.1. Since normality is not a transitive relation,
the first step is to add the conditions for normality of ${H_1}$:\\
 ${(m{x^a}{y^b}{z^c})^m}$, ${(r{x^d}{y^e}{z^f})^m}$, ${(m{x^a}{y^b}{z^c})^r}$, ${(r{x^d}{y^e}{z^f})^r}$,
 ${(m{x^a}{y^b}{z^c})^x}$, ${(r{x^d}{y^e}{z^f})^x}$, \\ ${(m{x^a}{y^b}{z^c})^y}$, ${(r{x^d}{y^e}{z^f})^y}$, 
${(m{x^a}{y^b}{z^c})^z}$,  ${(r{x^d}{y^e}{z^f})^z}$ $ \in {H_1} \cap T$. 

We obtain that  ${y^{2b}},{x^{2a}}{z^{2c}},{y^2}$, ${x^{2d}}{z^{2f}},{y^{2e}},$
${x^2},{z^2} \in {H_1} \cap T$. \\The conditions ${x^{2a}}{z^{2c}},{y^{2e}} \in {H_1} \cap T$ already being involved in the process, 
we may omit them while forming the second system of equations:  \\

$C_2 = \left\{ \begin{gathered}
  0 = {\alpha _5}g,2b = h{\alpha _5} + j{\beta _5},0 = i{\alpha _5} + k{\beta _5} + l{\gamma _5} \hfill \\
  0 = {\alpha _6}g,2 = h{\alpha _6} + j{\beta _6},0 = i{\alpha _6} + k{\beta _6} + l{\gamma _6}, \hfill \\
  2d = {\alpha _7}g,0 = h{\alpha _7} + j{\beta _7},2f = i{\alpha _7} + k{\beta _7} + l{\gamma _7} \hfill \\
  2 = {\alpha _8}g,0 = h{\alpha _8} + j{\beta _8},0 = i{\alpha _8} + k{\beta _8} + l{\gamma _8} \hfill \\
  0 = {\alpha _9}g,0 = h{\alpha _9} + j{\beta _9},2 = i{\alpha _9} + k{\beta _9} + l{\gamma _9} \hfill \\ 
\end{gathered}  \right\}$.\\

Solving the system that consists of equations given in $C'_1 $ and $ C_2$, we conclude that
the contribution to the normal zeta function coming from $ H_1 $ is:

 $ 1 + 8 \cdot {2^{ - s}} + 16 \cdot {2^{ - s}} + 4 \cdot {2^{ - s}} + 16 \cdot {4^{ - s}} + 32 \cdot {4^{ - s}} + 64 \cdot {4^{ - s}} + 64 \cdot {8^{ - s}} \\ 
 = 1 + 28 \cdot {2^{ - s}} + 112 \cdot {4^{ - s}} + 64 \cdot {8^{ - s}}. $

It is easily seen that a normal subgroup of ${G_{{2_1}}}$ is also a normal subgroup of $G=P2/m$. In this case, 
consideration of subgroups that are not contained in  ${G_{{2_2}}} $, ${G_{{2_3}}}$, ${G_3} $ 
and the facts from Theorem 1.1 yield the normal zeta function contribution:\\  ${2^{ - s}}(1 + 11 \cdot {2^{ - s}} + 12 \cdot {2^{ - 2s}})\zeta (s)\zeta_1 (s)$.

A normal subgroup of ${G_{{2_2}}\cong {P2}}$ is also a normal subgroup of  ${G=P2/m}$. 
Counting only those groups that are not contained in  ${G_{{2_3}}}$ and ${G_3}$ and using Theorem 1.1., 
we get the contribution to the normal zeta function of group $P2/m$:\\ 
${2^{ - s}}\left( {1 + 13 \cdot {2^{ - s}} + 22 \cdot {2^{ - 2s}} + 4 \cdot {2^{ - 3s}}} \right)\zeta (s)$.

The normal subgroup of ${G_{{2_3}}\cong P{\bar 1}}$ of the form 
$H_4=\left\langle {{mr}{x^a}{y^b}{z^c},{x^d}{y^e}{z^f},{y^g}{z^h},{z^i}} \right\rangle $, with  $0  \le a < d, 0 \le b,e < g , 0 \le c,f,h < i $,  is a 
normal subgroup of  ${G=P2/m}$. Again,  Theorem 1.1  and the groups that are not contained in  ${G_3} $
yield the share ${2^{ - s}}\left( {1 + 14 \cdot {2^{ - s}} + 28 \cdot {2^{ - 2s}} + 8 \cdot {2^{ - 3s}}} \right)$.

Denote a subgroup of  ${G_3} = \left\langle {x,y,z} \right\rangle$ by $H_5$. It takes the form 
${H_5} = \left\langle {{x^a}{y^b}{z^c},{y^d}{z^e},{z^f}} \right\rangle $, where we assume $0 < a,0 \le b < d,0 \le c,e < f$. 
Based on the conditions of normality, we deduce ${y^{2d}}, {y^{2b}} \in {H_5}$  and  another set of constraints:\\

$C = \left\{ \begin{gathered}
  0 = a{\alpha _1},2b = b{\alpha _1} + d{\beta _1},0 = c{\alpha _1} + e{\beta _1} + f{\gamma _1} \hfill \\
  0 = a{\alpha _2},2d = b{\alpha _2} + d{\beta _2},0 = c{\alpha _2} + e{\beta _2} + f{\gamma _2} \hfill \\ 
\end{gathered}  \right\}.$\\

The equations $0 = a{\alpha _1},2b = b{\alpha _1} + d{\beta _1}$ and the assumption that $d$ is even imply that $b$ can be $0$ or $\frac{d}{2}$. 
If we assume that $d$ is odd, then $b=0$ is the only choice for $b$. Similarly,
the equations $ 0 = a{\alpha _2},2d = b{\alpha _2} + d{\beta _2},0 = c{\alpha _2} + e{\beta _2} + f{\gamma _2}$ and an assumption that $f$ is odd
yield that there is only one choice for $e$. However, if we assume that  $f$ is even, then there exist two options for $e$. 
If $d$ and $f$ are both even and $b=0$, then there are two choices for $e$; if $b = \frac{d}{2}$, then  $e$ has to be 0.

We have the following contribution:\\

$\begin{gathered}
  {2^{ - 2s}}(3\sum\limits_{d,f \in 2\mathbb{N},a \in \mathbb{N}} {{a^{ - s}}{d^{ - s}}{f^{ - s}}f  }+ 2\sum\limits_{d \in 2\mathbb{N},f \in \mathbb{N}',a \in \mathbb{N}} {{a^{ - s}}{d^{ - s}}{f^{ - s}}f  }+\\+ 2\sum\limits_{d \in \mathbb{N}',f \in 2\mathbb{N},a \in \mathbb{N}} {{a^{ - s}}{d^{ - s}}{f^{ - s}}f}  +  \sum\limits_{d,f \in \mathbb{N}',a \in \mathbb{N}} {{a^{ - s}}{d^{ - s}}{f^{ - s}}f} ) \\= {2^{ - 3s}}(3 + {2^s}){\zeta ^2}(s)\zeta_1 (s). \hfill \\ 
\end{gathered} $\\

Finally, we obtain the explicit expression for the normal zeta function of group ${P_{2/m}}$:  \\

$  \zeta _{{P_{2/m}}}^ \triangleleft (s) = 1 + 28 \cdot {2^{ - s}} + 112 \cdot {4^{ - s}} + 64 \cdot {8^{ - s}} + {2^{ - s}}(1 + 11 \cdot {2^{ - s}} +12 \cdot {2^{ - 2s}}) \cdot \\
\zeta (s)\zeta_1 (s)+ 
2^{ - s}\cdot (1 + 13 \cdot {2^{ - s}} + 22 \cdot {2^{ - 2s}} + 4 \cdot {2^{ - 3s}})\zeta (s)  \\+ {2^{ - s}}\left( {1 + 14 \cdot {2^{ - s}} + 28 \cdot {2^{ - 2s}} + 8 \cdot {2^{ - 3s}}} \right) +  {2^{ - 3s}}(3 + {2^s}){\zeta ^2}(s)\zeta_1 (s)\\ = 1 + 29 \cdot {2^{ - s}} + 126 \cdot {4^{ - s}} + 92 \cdot {8^{ - s}} + 8 \cdot {16^{ - s}} \\
  + {2^{ - s}}\left( {1 + 13 \cdot {2^{ - s}} + 22 \cdot {2^{ - 2s}} + 4 \cdot {2^{ - 3s}}} \right)\zeta (s) + {2^{ - s}}(1 + 11 \cdot {2^{ - s}} + 12 \cdot {2^{ - 2s}})\zeta (s)\zeta_1 (s)\\ 
   + {2^{ - 3s}}(3 + {2^s}){\zeta ^2}(s)\zeta_1 (s). \hfill \\ $

\section{\bf Subgroups of finite index in $P2/m$}
\mbox{}

In the sequel,  $d(n)$ denotes the number of all positive divisors of a positive integer $n$ and $\sigma (n)$ is the sum of all positive divisors for a 
positive integer $n$, i. e. $\sigma (n) = \sum\limits_{\left. l \right|n} l ,$ as usual. The answer we are looking 
for is contained in the next two propositions. Their validity readily follows from Theorem 2.1. and Theorem 3.1.
and the fact that the product of two Dirichlet series 
$\sum\limits_{n \in \mathbb {N}}{f(n)}{n^{ - s}}$ and $\sum\limits_{n \in \mathbb {N}}{g(n)}{n^{ - s}}$, where $f$ and $g$ are two arithmetic functions, is a
Dirichlet series $\sum\limits_{n \in \mathbb {N}}{h(n)}{n^{ - s}}$  with the coefficients $h(n) = (f * g)(n) = \sum\limits_{\left. l \right|n} {f(l)g\left( {\frac{n}{l}} \right)}  = \sum\limits_{ab = n} {f(a)g(b)} $. The last sum extends over all positive divisors $l$ of $n$, or equivalently over all distinct pairs $(a, b)$ of positive integers whose product is $n$. Let us remember that $d(n)$ and $\sigma (n)$ are coefficients of Dirichlet series $\zeta (s)\zeta(s)$ and  $\zeta (s)\zeta_1(s)$ respectively.

\begin{prop} The number $a_n$ of all  subgroups of  index $n$ in the group $P2/m$ is given by the following expressions
\begin{enumerate}       

  \item  if $n$ is even, \[{a_n} = \left\{ \begin{array}{l}
\left. {n\sum\limits_{\left. l \right|n} {\sigma (l)}  + 10n\sum\limits_{\left. l \right|\left( {\frac{n}{2}} \right)} {\sigma (l)}  + \sum\limits_{\left. l \right|\left( {\frac{n}{2}} \right)} {l \cdot d(l)}  + \left( {\frac{n}{2} + 1} \right)\sum\limits_{\left. l \right|\left( {\frac{n}{2}} \right)} {l \cdot \sigma (l)} } \right\}\\
\,(n \equiv 2\, \vee n \equiv 6)\,(\bmod 8)\\
\\
\left. \begin{array}{l}
n\sum\limits_{\left. l \right|n} {\sigma (l)}  + 10n\sum\limits_{\left. l \right|\left( {\frac{n}{2}} \right)} {\sigma (l)}  + 9n\sum\limits_{\left. l \right|\left( {\frac{n}{4}} \right)} {\sigma (l)}  + \sum\limits_{\left. l \right|\left( {\frac{n}{2}} \right)} {ld(l)}  + 9\sum\limits_{\left. l \right|\left( {\frac{n}{4}} \right)} {ld(l)}  + \\
 + 8\sum\limits_{\left. l \right|\left( {\frac{n}{4}} \right)} {l\sigma (l) + } \left( {\frac{n}{2} + 1} \right)\sum\limits_{\left. l \right|\left( {\frac{n}{2}} \right)} {l\sigma (l)} 
\end{array} \right\}\\
\,\,n \equiv 4\,(\bmod 8)\\
\\
\left. \begin{array}{l}
n\sum\limits_{\left. l \right|n} {\sigma (l)}  + 10n\sum\limits_{\left. l \right|\left( {\frac{n}{2}} \right)} {\sigma (l)}  + 9n\sum\limits_{\left. l \right|\left( {\frac{n}{4}} \right)} {\sigma (l)}  + \sum\limits_{\left. l \right|\left( {\frac{n}{2}} \right)} {ld(l)}  + 9\sum\limits_{\left. l \right|\left( {\frac{n}{4}} \right)} {ld(l)}  + \\
 + 8\sum\limits_{\left. l \right|\left( {\frac{n}{4}} \right)} {l\sigma (l) + } \left( {\frac{n}{2} + 1} \right)\sum\limits_{\left. l \right|\left( {\frac{n}{2}} \right)} {l\sigma (l)}  + 6\sum\limits_{\left. l \right|\left( {\frac{n}{8}} \right)} {ld(l)} 
\end{array} \right\}\,\\
n \equiv 0(\bmod 8)
\end{array} \right.\]

\item	if $n$ is odd,\[{a_n} = n\sum\limits_{\left. l \right|n} {\sigma (l)} \]
\end{enumerate}
In particular, ${a_p} = {p^2} + 2p$ for every odd prime $p$. 
\end{prop}

\begin {prop} The number $c_n$ of all  normal subgroups of  index $n$ in group  $P2/m$ reads:
\begin{enumerate}
\item ${c_1} = 1$
\item  if $n$ is odd and $n \ne 1,{c_n} = 0$
\item  if $n$ is even,\[{c_n} = \left\{ \begin{array}{l}
31,\,n = 2\\
155,n = 4\\
187,n = 8,\\
199,n = 16\\
40 + \sigma (n/2) + 11 \cdot \sigma (n/4) + 12 \cdot \sigma (n/8) + 3 \cdot \sum\limits_{l\left| {\left( {\frac{n}{8}} \right)} \right.} {\sigma (l)} \, + \sum\limits_{l\left| {\left( {\frac{n}{4}} \right)} \right.} {\sigma (l)} \,,\,\,\,\,\\
\left( {\,n \equiv 0(\bmod 16) \wedge n \ne 16} \right)\\
\\
36 + \sigma (n/2) + 11 \cdot \sigma (n/4) + 12 \cdot \sigma (n/8) + 3 \cdot \sum\limits_{l\left| {\left( {\frac{n}{8}} \right)} \right.} {\sigma (l)} \, + \sum\limits_{l\left| {\left( {\frac{n}{4}} \right)} \right.} {\sigma (l)} ,\\
\left( {\,n \equiv 8(\bmod 16) \wedge n \ne 8} \right)\\
\\
14 + \sigma (n/2) + 11 \cdot \sigma (n/4) + \sum\limits_{l\left| {\left( {\frac{n}{4}} \right)} \right.} {\sigma (l)} ,\\
\left( {\,(n \equiv 4 \vee n \equiv 12)(\bmod 16) \wedge n \ne 4} \right)\\
1 + \sigma (n/2),\left( {\,(n \equiv 2 \vee n \equiv 6 \vee n \equiv 10 \vee n \equiv 14)(\bmod 16) \wedge n \ne 2} \right)
\end{array} \right.\]

\end{enumerate}

\end {prop}

\mbox{}

\begin{prop} $P2/m$ is a group of degree 3.

\end {prop}

{\bf Proof.} Recall that the degree of a group $G$ is defined by $\deg (G) = \lim \sup \frac{{\log {a_n}(G)}}{{\log n}}$, 
where ${a_n}(G)$ is the number of subgroups of index $n$ in $G$. In other words, deg($G$) is the ``smallest'' positive
 real number $c$ such ${a_n}(G) = O({n^{c + \varepsilon }})$, for all $\varepsilon  > 0$ and all $n$. 

 Proposition 4.1. implies  ${a_n} = O(n\sum\limits_{\left. l \right|n} {l\sigma (l)} )$.
By Robin's inequality, we have $\sigma (l) = O(l\log\log l)$. Hence,\\
$\sum\limits_{\left. l \right|n} {l\sigma (l)} = O(\log\log n \sum\limits_{\left. l \right|n} {l^2})
= O(n\log\log n \sigma (n)) = O(n^2(\log\log n)^2)$.\\
Thus, ${a_n} = O(n^3(\log\log n)^2) = O({n^{3 + \varepsilon }})$ for every $\varepsilon  > 0$ and every $n$.

On the other hand, let us have a look at the subsequence ${a_{2p}}$, where $p$ runs through the prime numbers. 
Obviously, $2p \equiv 2\, \vee 2p \equiv 6\,(\bmod 8)$. Now, Proposition 4.1. and straightforward calculations 
yield $a_{2p} = p^3 + 30p^2 + 60p + 2$. Then \\
$\lim \sup \frac{{\log {a_{2p}}}}{{\log {2p}}} = 3$. Thus, $\deg (P2/m) = 3$.\\

\section{\bf Subgroup growth in $P2/m$}

\mbox{}

\begin {thm} $\sum\limits_{n \le x} {{a_n}}  = \frac{{{x^4}{\pi ^2}}}{{384}}\zeta (3) + O({x^3}\ln x)$.
\end{thm}

In the proof of the theorem, we shall make use of the following lemma.\\

\begin {lem} $\sum\limits_{n \le x} {\sum\limits_{\left. q \right|n} {q\sigma (q)} }  = \frac{{{\pi ^2}}}{{18}}\zeta (3){x^3} + O({x^2}\ln x)$.
\end{lem}

{\bf Proof of Lemma.} Note that

$ \sum\limits_{nmh \le x} {n{m^2}}  = \sum\limits_{n \le x} {\sum\limits_{\left. q \right|n} {q\sigma (q)} }  
= \sum\limits_{d \le x} {\sum\limits_{q \le \frac{x}{d}} {q\sigma (q)}}.$ \\

By Abel's summation formula, the right-hand side further equals

 $\sum\limits_{d \le x} {\left\{ {\frac{x}{d}\left( {\sum\limits_{q \le \frac{x}{d}} {\sigma (q)} } \right) - 1 \cdot 1 - \int\limits_1^{\frac{x}{d}} {\sum\limits_{q \le t} {\sigma (q)dt}} } \right\}}.$\\

Using the well known fact 

${\sum\limits_{q \le t} {\sigma (q)}}  = \frac{{{\pi ^2}}}{{2}}\zeta (2){t^2} + O({t}\ln t)$,

we transform the above expression into

$\sum\limits_{d \le x} {\left\{ {\frac{x}{d}\left( {\frac{1}{2}\zeta (2){{\left( {\frac{x}{d}} \right)}^2} + O\left( {\frac{x}{d}\ln \left( {\frac{x}{d}} \right)} \right)} \right) - 1 - \int\limits_1^{\frac{x}{d}} {\left( {\frac{1}{2}\zeta (2){t^2} + O\left( {t\ln t} \right)} \right)dt} } \right\}}  =$ 

$=\frac{{{x^3}\zeta (2)}}{2}\sum\limits_{d \le x} {\frac{1}{{{d^3}}}}  + O({x^2}\ln x) - \sum\limits_{d \le x} {\left. {\left( {\frac{1}{2}\zeta (2)\frac{{{t^3}}}{3} + O(\frac{{{t^2}}}{2}\ln t - \frac{{{t^2}}}{4})} \right)} \right|_1^{\frac{x}{d}}} $

$ = \frac{{{x^3}\zeta (2)}}{2}\sum\limits_{d \le x} {\frac{1}{{{d^3}}}}  + O({x^2}\ln x)
 - \sum\limits_{d \le x} {\left( {\frac{1}{2}}\zeta (2)\frac{{{x^3}}}{{3{d^3}}}+ O(\frac{{{x^2}}}{{2{d^2}}}\ln \frac{x}{d}) \right)} $

$= \frac{1}{3}\zeta (2)\zeta (3){x^3} + O({x^2}\ln x) =\frac{{{\pi ^2}}}{{18}}\zeta (3){x^3} + O({x^2}\ln x)$.\\

{\bf {Proof of Theorem 5.1.}} Accordig to Proposition 4.1, one has

$\sum\limits_{n \le x} {{a_n}}  = \sum\limits_{{n_1}{n_2}{n_3} \le x} {{n_1}{n_2}n_3^2}  + 20 \cdot \sum\limits_{2{n_1}{n_2}{n_3} \le x} {{n_1}{n_2}n_3^2}  + 36 \cdot \sum\limits_{4{n_1}{n_2}{n_3} \le x} {{n_1}{n_2}n_3^2}  +$

$+\sum\limits_{2{n_1}{n_2}{n_3} \le x} {{n_1}{n_2}\,}  + 9 \cdot \sum\limits_{4{n_1}{n_2}{n_3} \le x} {{n_1}{n_2}\,} 
 + 6 \cdot \sum\limits_{8{n_1}{n_2}{n_3} \le x} {{n_1}{n_2}\,}  + \sum\limits_{2{n_1}{n_2}{n_3} \le x} {{n_1}n_2^2\,}  +$\\

$+ 8 \cdot \sum\limits_{4{n_1}{n_2}{n_3} \le x} {{n_1}n_2^2\,}  + \sum\limits_{2{n_1}{n_2}{n_3} \le x} {{n_1}n_2^2n_3^3}.\\
$

The leading term is $\sum\limits_{2{n_1}{n_2}{n_3} \le x} {{n_1}n_2^2n_3^3} $. By Abel's partial summation formula,

$\sum\limits_{{n_1}{n_2}{n_3} \le x} {{n_1}n_2^2n_3^3} = \sum\limits_{n \le x} {n\sum\limits_{\left. q \right|n} {q\sigma (q)} }= x \cdot \sum\limits_{n \le x} {\sum\limits_{\left. q \right|n} {q\sigma (q)} }  - 1 \cdot 1 - \int\limits_1^x {\sum\limits_{n \le t}} {\sum\limits_{\left. q \right|n} {q\sigma (q)} } dt$.\\

Lemma implies that this is equal to 

$\frac{{{\pi ^2}}}{{18}}\zeta (3){x^4}
+O({x^3}\ln x)- \int\limits_1^x {\left( {\frac{{{\pi ^2}}}{{18}}\zeta (3){t^3} + O({t^2}\ln t)} \right)} dt =$

$\begin{array}{l}
 = \frac{{{\pi ^2}}}{{18}}\zeta (3){x^4} + O({x^3}\ln x) - \left. {(\frac{{{\pi ^2}}}{{18}}\zeta (3)\frac{{{t^4}}}{4} + O(\frac{{{t^3}}}{3}\ln t - \frac{{{t^3}}}{9}))} \right|_1^x = \\
 = \frac{{{\pi ^2}}}{{24}}\zeta (3){x^4} + O({x^3}\ln x).
\end{array}$\\

Hence,  $\sum\limits_{2{n_1}{n_2}{n_3} \le x} {{n_1}n_2^2n_3^3} =\sum\limits_{{n_1}{n_2}{n_3} \le \frac{x}{2}} {{n_1}n_2^2n_3^3} = \frac{{{x^4}{\pi ^2}}}{{384}}\zeta (3) + O({x^3}\ln x)$.

\begin {thm} $\sum\limits_{n \le x} {{c_n}}  = \left( {\frac{3}{{32}} + \frac{{7{\pi ^2}}}{{4608}}} \right){x^2}{\pi ^2} + O(x{\ln ^2}x)$
\end {thm}
{\bf Proof.} By Proposition 4.2, we have

$\sum\limits_{n \le x} {{c_n}}  = 256 + \sum\limits_{2{n_1}{n_2} \le x} {{n_1}}  + 11 \cdot \sum\limits_{4{n_1}{n_2} \le x} {{n_1}}  + 12 \cdot \sum\limits_{8{n_1}{n_2} \le x} {{n_1}}  + 3 \cdot \sum\limits_{8{n_1}{n_2}{n_3} \le x} {{n_1}\,}+$

$  + \sum\limits_{4{n_1}{n_2}{n_3} \le x} {{n_1}\, + } \sum\limits_{2{n_1}{n_2} \le x} 1  
 + 13 \cdot \sum\limits_{4{n_1}{n_2} \le x} 1  + 22\sum\limits_{8{n_1}{n_2} \le x} {1 + 4 \cdot \sum\limits_{16{n_1}{n_2} \le x} 1 }\\ 
$

Recall once again that $\sum\limits_{nm \le x} n  = \frac{{\zeta (2)}}{2} \cdot {x^2} + O(x\ln x)$.

On the other hand,

$\sum\limits_{nmh \le x} {n = } \frac{{{\pi ^4}}}{{72}}{x^2} + O(x{\ln ^2}x).$

Indeed,

 $\sum\limits_{nmh \le x} n  = \sum\limits_{n \le x} {\sum\limits_{\left. q \right|n} {\sigma (q)} }  
= \sum\limits_{d \le x} {\sum\limits_{q \le \frac{x}{d}} {\sigma (q)} } 
= \sum\limits_{d \le x} {\left\{ {\frac{1}{2}\zeta (2){{\left( {\frac{x}{d}} \right)}^2} + O\left( {\frac{x}{d}\ln \left( {\frac{x}{d}} \right)} \right)} \right\}} $

$ = \frac{1}{2}\zeta (2){x^2}\sum\limits_{d \le x} {\frac{1}{{{d^2}}} + } O(x\ln x\sum\limits_{d \le x} {\frac{1}{d} - x} \sum\limits_{d \le x} {\frac{{\ln d}}{d}} )$

$= \frac{1}{2}{\zeta ^2}(2){x^2} + O(x{\ln ^2}x) = \frac{{{\pi ^4}{x^2}}}{{72}} + O(x{\ln ^2}x)$.

Thus, we get

$\sum\limits_{n \le x} {{c_n}}  = \frac{1}{2}\zeta (2){\left( {\frac{x}{2}} \right)^2} + O\left( {\frac{x}{2}\ln \left( {\frac{x}{2}} \right)} \right) + 11 \cdot \left( {\frac{1}{2}\zeta (2){{\left( {\frac{x}{4}} \right)}^2} + O\left( {\frac{x}{4}\ln \left( {\frac{x}{4}} \right)} \right)} \right) + $

$ + 12 \cdot \left( {\frac{1}{2}\zeta (2){{\left( {\frac{x}{8}} \right)}^2} + O\left( {\frac{x}{8}\ln \left( {\frac{x}{8}} \right)} \right)} \right) + 3 \cdot \left( {\frac{{{\pi ^4}}}{{72}}{{\left( {\frac{x}{8}} \right)}^2} + O\left( {\frac{x}{8}{{\ln }^2}\left( {\frac{x}{8}} \right)} \right)} \right) +$ 

$ + \left( {\frac{{{\pi ^4}}}{{72}}{{\left( {\frac{x}{4}} \right)}^2} + O\left( {\frac{x}{4}{{\ln }^2}\left( {\frac{x}{4}} \right)} \right)} \right) = \frac{{3{x^2}{\pi ^2}}}{{32}} + \frac{{7{x^2}{\pi ^4}}}{{4608}} + O(x{\ln ^2}x).
$\\

\end{document}